\documentclass[10pt,a4paper,oneside,leqno]{article}
\usepackage{amsfonts,amssymb}
\usepackage{amscd}
\newtheorem{defn}{Definition}
\newtheorem{conven}{Convention}

\newtheorem{observen}{Observation}

\vspace{0.5cm}
 4
\font\ebf=cmbx8
\font\erm=cmr8

\usepackage{graphicx, graphics}
\usepackage{color}
\newcommand{\os}{\oplus\!\!\to}
\newcommand{\cs}{\copyright\!\!\to}

\begin{document}
\begin{center}
	\noindent { \textsc{ Cobweb Posets and  KoDAG  Digraphs are Representing Natural Join of Relations, their di-Bigraphs and the Corresponding  Adjacency Matrices.}}  \\ 
	\vspace{0.3cm}
	\vspace{0.3cm}
	\noindent Andrzej Krzysztof Kwa\'sniewski \\
	\vspace{0.2cm}
	\noindent {\erm Member of the Institute of Combinatorics and its Applications  }\\
{\erm High School of Mathematics and Applied Informatics} \\
	{\erm  Kamienna 17, PL-15-021 Bia\l ystok, Poland }\noindent\\
	\noindent {\erm e-mail: kwandr@gmail.com}\\
	\vspace{0.4cm}
\end{center}

\noindent {\ebf Abstract:}
\vspace{0.1cm}
\noindent {\small Natural join of di-bigraphs (directed bi-parted graphs) and their corresponding  adjacency matrices is defined and then applied to investigate the so called cobweb posets and their $Hasse$ digraphs called $KoDAGs$. $KoDAGs$  are special \textbf{o}rderable \textbf{D}irected \textbf{A}cyclic \textbf{G}raphs which are cover relation digraphs of cobweb posets introduced by the author few years ago. $KoDAGs$   appear to be distinguished  family of $Ferrers$ digraphs which are  natural join of  a corresponding ordering chain of one direction directed cliques called di-bicliques. These digraphs  serve to represent faithfully  corresponding relations of  arbitrary arity so that all relations of arbitrary arity are their subrelations. Being this $chain -way$  complete (compare with \textbf{K}ompletne , \textbf{K}uratowski  $K_{n,m}$ bipartite  graphs)  their DAG denotation is accompanied  with the letter \textbf{K }in front of  descriptive abbreviation oDAG. 
\vspace{0.1cm}
\noindent The way to join bipartite digraphs of binary into multi-ary relations is the natural join operation either on relations or their digraph representatives. This natural join operation  is denoted here by $\os$ symbol deliberately referring - in a reminiscent manner - to the direct sum  $\oplus$ of adjacency  matrices as it becomes the case for disjoint di-bigraphs.}

\vspace{0.3cm}

\noindent Key Words: posets, graded digraphs, Ferrers dimension,  natural join

\vspace{0.1cm}

\noindent AMS Classification Numbers: 06A06 ,05B20, 05C7  

\vspace{0.1cm}

\noindent  affiliated to The Internet Gian-Carlo Polish Seminar:

\noindent \emph{http://ii.uwb.edu.pl/akk/sem/sem\_rota.htm}

%%%%%%%%%%%%%%%%%%%%%%%%%%%%%%%%%%%%%%%%%%%%%%%%%%%%%%%%%%%%%%%%%%
\section{Introduction to coweb posets}

\subsection{Notation}
One may identify and interpret some classes  of digraphs in terms of their associated posets.  (see \cite{1} Interpretations in terms of posets Section 9)

\begin{defn}[see \cite{1}]
Let $D = (\Phi,\prec)$      be a digraph.   $w,v \in \Phi$ are said to be equivalent iff  there exists a directed path containing both $w$ and $v$ vertices.  We then write:  $v \sim w$  for such pairs and  denote by $[v]$  the  $\sim$ equivalence class of  $v \in \Phi$.
\end{defn}

\begin{defn}[see \cite{1}]
The poset $P(D)$ associated to $D = (\Phi,\prec)$    is the poset $P(D)= (\Phi / \sim , \leq)$    where
 $[v] \leq [w]$   iff  there exists a directed path from a vertex  $x \in [v]$  to a vertex  $y \in [w]$.
\end{defn}

\vspace{0.2cm}
\noindent \textbf{The graded  digraphs case:}

\vspace{0.2cm}
\noindent \textbf{If  $D = (\Phi,\prec )$     is graded digraph then $D = (\Phi, \prec )$}   is necessarily \textbf{acyclic}.
Then no two elements of  $D = (\Phi,\prec  )$ are   $\sim$ equivalent and thereby
$P(D) = (\Phi / \sim , \leq)$    associated to  $D = (V,\prec  )$    \textbf{is equivalent to}:    $P(D) \equiv (\Phi  , \leq)$   =  transitive, reflexive  closure  of $D = (\Phi,\prec  )$.

\vspace{0.2cm}
\noindent The cobweb posets where introduced in several paper (see \cite{2}-\cite{6} and references therein) in terms of their poset [Hasse] diagrams.  Here we deliver their equivalent  definition preceded by preliminary notation and nomenclature.

\vspace{0.2cm}
\noindent \textbf{Notation :  nomenclature, di-bicliques  and natural join}

\vspace{0.2cm}
\noindent In order to proceed proficiently we  adopt the following.

\begin{defn}
A digraph $D = (\Phi,\prec\!\!\cdot)$ is  transitive irreducible  iff  transitive $reduction(D)  = D$.
\end{defn}

\begin{defn}
A poset $P(D) = (\Phi, \leq)$ is associated to  a graded digraph $D = (\Phi,\prec  )$   iff   $P(D)$ is the  transitive, reflexive  closure  of $D = (\Phi, \prec )$ .
\end{defn}

\vspace{0.2cm}
\noindent \textbf{Obvious}.

\noindent $D = (\Phi,\prec\!\!\cdot)$ is  transitive irreducible  iff  transitive $reduction(D)  = D$  iff  $D = (\Phi,\prec\!\!\cdot  )$   is Hasse diagram of the poset $P(D) = (\Phi, \leq)$   associated to  $D  \equiv  D = (\Phi,\prec\!\!\cdot  )$   is cover relation $\prec\!\!\cdot$  digraph $\equiv$    $D = (\Phi,\prec\!\!\cdot  )$  is  $P(D) = (\Phi, \leq)$    poset diagram.

%%%%%%%%%%%
\subsection{ Further on we adopt also the following nomenclature.}

We shall use until stated otherwise the convention: $N = \{1,2,...,k,...\}$ .   $n \in N \cup \{\infty\}$.  
The Cartesian product   $\Phi_1\times...\times\Phi_k$  of pairwise disjoint sets $\Phi_1, ... , \Phi_k$   is a  $k$-ary relation, called sometimes the universal relation and here now on  \textbf{K}ompletna relation or \textbf{K}-relation, (in Professor  \textbf{K}azimierz \textbf{K}uratowski native language this means complete). The purpose of introducing the letter $K$  is to distinguish  in what follows [ for $k = 2$  ] from  complete digraphs notions  established content.

\begin{conven}[identification]
The binary relation $E \subseteq X\times Y$  is being here identified with its bipartite digraph representation    $B = (X \cup Y, E)$.
\end{conven}

\noindent \textbf{Notation}  $\stackrel{\rightarrow}{K_{m,n}}\equiv B = (X \cup Y, E)$ if  $|X|= m$ ,   $|Y| = n$. Colligate with \textbf{K}uratowski and $K_{m,n}$.

\vspace{0.4cm}
\noindent \textbf{Comment 1.}

\noindent Complete  $n$-vertex \textbf{graphs} for which  all pairs of vertices are adjacent  are denoted by $K_n$, The letter $K$ had been chosen  in honor of Professor \textbf{K}azimierz  \textbf{K}uratowski, a distinguished pioneer in graph theory. The corresponding two widely used  concepts for digraphs are called  complete  digraphs or complete symmetric digraph  in which every two different vertices are joined by an arc  and complete oriented graphs  i.e. tournament graphs.

\vspace{0.2cm}
\noindent The binary $K$-relation  $E  = X\times Y$  equivalent  to  bipartite digraph  $B = ( X \cup Y, E) \equiv \stackrel{\rightarrow}{K_{m,n}}$  is called  from now on a \textbf{di-biclique}  following  \cite{6}.

\vspace{0.4cm}
\noindent \textbf{Example of}  di-bicliques obtained from  \textbf{bicliques} :  See Fig. 1. 

\vspace{0.1cm}
\noindent  If you  imagine arrows $\to$ left to the right - you would see two examples  of  \textbf{di-bicliques}         
 
%%%%%%% image
\begin{figure}[ht]
\begin{center}
	\includegraphics[width=50mm]{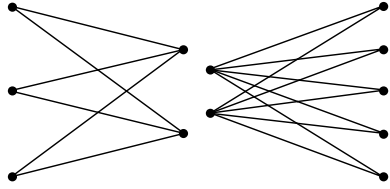}
	\caption{Examples of di-bicliques if edges are replaced by arrows  of join direction \label{fig:1}}
\end{center}
\end{figure}

\noindent if you  imagine arrows $\leftarrow$ right to the left, you would see another examples of  \textbf{di-bicliques}.         

\begin{conven}[recall]
The binary relation $E \subseteq X\times Y$  is identified with its  bipartite digraph  $B = ( X \cup Y, E)$ unless
otherwise denoted distinctively  deliberately.
\end{conven}

\vspace{0.2cm}
\noindent \textbf{The natural join.}

\vspace{0.2cm}
\noindent The natural join operation is a binary operation  like $\Theta$  \textbf{operator in computer science} denoted here by $\os$ symbol  deliberately referring - in a quite  reminiscent manner - to direct sum  $\oplus$  of adjacency Boolean  matrices and - as matter of fact and  in effect - to direct the sum  $\oplus$  of corresponding biadjacency [reduced]  matrices of digraphs under natural join.

\noindent  $\os$ is a natural operator for sequences construction .  $\os$ operates on multi-ary relations according to the scheme:   $(n+k)_{ary}  \os (k+m)_{ary}$  =   $(n+ k +m)_{ary}$

\noindent For example:  $(1+1)_{ary} \os(1+1)_{ary} = (1+ 1 +1)_{ary}$ , binary $\os$ binary = ternary.

\vspace{0.2cm} 
\noindent Accordingly an action of  $\os$  on these multi-ary relations' digraphs adjacency matrices is to be designed soon in what follows.

\vspace{0.2cm}
\noindent \textbf{Domain-Codomain $F$-sequence condition} $\mathrm{dom}(R_{k+1}) = \mathrm{ran} (R_k)$,  $k=0,1,2,...$ .  

\vspace{0.4cm}
\noindent Consider any  natural number valued sequence $F = \{F_n\}_{n\geq 0}$. Consider then any  chain of binary relations  defined on pairwise disjoint finite sets with cardinalities appointed by  $F$ -sequence elements values. For that to start we specify at first   a relations'  domain-co-domain $F$ - sequence.   

\vspace{0.4cm}
\noindent \textbf{Domain-Codomain  $F$-sequence $(|\Phi_n| = F_n )$}

$$
	\Phi_0,\Phi_1,...\Phi_i,...\ \ \Phi_k\cap\Phi_n = \emptyset \ \ for \ \ k \neq n, |\Phi_n|=F_n; \ \ i,k,n=0,1,2,...
$$

\noindent Let $\Phi=\bigcup_{k=0}^n\Phi_k$  be the corresponding ordered partition  [ anticipating -  $\Phi$ is the vertex set of  $D = (\Phi,\prec\!\!\cdot$ )   and  its transitive, reflexive closure $(\Phi, \leq)$] .     Impose $\mathrm{dom} (R_{k+1}) = \mathrm{ran} (R_k)$ condition , $k\in N \cup \{\infty\}$. What we get is binary relations chain.

\begin{defn} [Relation`s chain]
Let   $\Phi=\bigcup_{k=0}^n\Phi_k$ ,   $\Phi_k  \cap \Phi_n = \emptyset$  for  $k  \neq  n$ be  the ordered partition of the set  $\Phi$ .
\end{defn}

\vspace{0.2cm}
\noindent Let a sequence of binary relations be given such that 
$$
	R_0,R_1,...,R_i,...,R_{i+n},...,\ \ R_k\subseteq\Phi_k\times\Phi_{k+1},\ \ \mathrm{dom}(R_{k+1}) = \mathrm{ran}(R_k).
$$

\noindent Then the sequence  $\langle R_k\rangle_{k\geq 0}$    is called natural join  (binary) \textbf{relation's chain}. 
Extension to varying arity relations' natural join chains is straightforward.

\noindent As necessarily $\mathrm{dom}(R_{k+1}) = \mathrm{ran}(R_k)$  for relations' natural join chain any given  binary relation's chain is not just a sequence  therefore we use "link to link " notation for $k, i , n = 1,2,3,...$ ready for  relational data basis applications: 
$$
 R_0 \os R_1 \os ... \os R_i \os ... \os R_{i+n},... is\ an\  F-chain\ of\  binary\ relations
$$

\noindent where  $\os$  denotes natural join of relations as well as both  natural join of their bipartite digraphs and the natural join of their representative adjacency matrices (see the Section 3.).

\vspace{0.2cm}
\noindent Relation's $F$-chain  naturally  represented by [identified with] the chain of theirs  \textbf{bipartite digraphs}

$$
	{ R_0 \os R_1 \os ... \os R_i \os ... \os R_{i+n},... \Leftrightarrow 
	\atop
	\Leftrightarrow  B_0 \os B_1 \os ... \os B_i \os ... \os B_{i+n},...
	}
$$

\vspace{0.2cm}
\noindent results in \textbf{$F$-partial ordered set} $\langle\Phi,\leq\rangle$    with its Hasse digraph representation  looking like 
specific "cobweb"   image  [see figures below].

%%%%%%%%%%%%%%%
\subsection{ Partial order  $\leq$}

The partial order relation $\leq$ in the set of all points-vertices is determined  uniquely by the  above equivalent $F$- chains.  Let  $x,y \in  \Phi=\bigcup_{k=0}^n\Phi_k$ and let  $k, i  = 0,1,2,...$.   Then 

\begin{equation}\label{eq:leq}
	x\leq y \Leftrightarrow \forall_{x\in\Phi} : x\leq x \vee \Phi_i\ni x < y \in \Phi_{i+k}\ iff\ x(R_i\copyright...\copyright R_{i+k-1})y
\end{equation}

\noindent where  "$\copyright$"  stays for [Boolean] composition of binary relations.

\vspace{0.4cm}
\noindent \textbf{Relation  ($\leq$) defined equivalently }: 

\vspace{0.2cm}
\noindent 
$ x  \leq  y$  in $(\Phi,\leq)$  iff either  $x=y$ or  there exist a directed path from $x$ to $y; x,y \in \Phi$.

\vspace{0.4cm}
\noindent Let now $R_k = \Phi_k\times\Phi_{k+1},  k \in N \cup\{0\}$. For "historical" reasons \cite{2}-\cite{6}  we shall call such partial ordered set $\Pi = \langle\Phi,\leq\rangle$  the \textbf{cobweb poset} as  theirs  Hasse digraph representation  looks like specific "cobweb"   image  ( imagine and/or draw also their transitive and reflexive cover  digraph $\langle\Phi,\leq\rangle$. Cobweb? Super-cobweb ! ...- with fog droplets loops ?) .  

%%%%%%%%%%%%
\subsection{ Cobweb posets  ($\Pi = \langle\Phi,\leq\rangle$) }

\begin{conven}[recall]
The binary relation $E  \subseteq X\times Y$  is identified with its  bipartite digraph  $B = ( X \cup Y, E)\equiv  \stackrel{\rightarrow}{K_{m,n}}$ where  $|X|= m, |Y| = n$.
\end{conven}

\begin{defn} [cobweb poset]
Let $D = (\Phi, \prec\!\!\cdot )$  be a  transitive irreducible  digraph. Let  $n \in N \cup \{\infty\}$.  Let $D$ be a natural join $D = \os_{k=0}^n B_k$   of   di-bicliques  $B_k = (\Phi_k \cup \Phi_{k+1}, \Phi_k\times\Phi_{k+1} ) , n \in N \cup \{\infty\}$.    Hence the digraph  $D = (\Phi,\prec\!\!\cdot  )$    is graded. The poset $\Pi (D)$ associated to this graded digraph $D = (\Phi,\prec\!\!\cdot )$    is called a cobweb poset.
\end{defn}

\vspace{0.2cm}

\begin{conven}
In a case we want to underline that we deal with finite cobweb poset ( a subposet of appropriate - for example infinite $F$-cobweb poset $\Pi (D)$ ) we shall use a subscript and write  $P_n$ .
\end{conven}
\noindent See: \cite{2}-\cite{6}, \cite{10}, \cite{13}, \cite{18}.

\vspace{0.2cm}

\noindent\textbf{Comment 2.}
\vspace{0.2cm}

\noindent \textbf{Graded graph}  is a \textbf{natural join} of bipartite graphs that form a chain of consecutive levels  [i.e. graded \textbf{graphs'}  antichains]

\vspace{0.2cm}
\noindent \textbf{Graded digraph}  is a \textbf{natural join} of bipartite digraphs that form a chain of consecutive levels [i.e. graded \textbf{digraphs'} antichains]

\vspace{0.4cm}
\noindent \textbf{Comment 3. } ({\it Definition 6.  Recapitulation in brief.})

\noindent Cobweb poset is the poset $\Pi = \langle\Phi,\leq\rangle$, where  $\Phi = \bigcup_{k=0}^n$ and $\prec\!\!\cdot = \os_{k=0}^{n-1} \Phi_k\times\Phi_{k+1},  n \in N \cup \{\infty\}$.
\noindent Cobweb poset is the poset $\Pi = \langle\Phi,\leq\rangle$, where  $\Phi = \bigcup_{k=0}^n$  and $\prec\!\!\cdot = \os_{k=0}^{n-1} \stackrel{\rightarrow}{K_{k,k+1}},  n \in N \cup \{\infty\}$,  
\noindent where  $\leq$ is the  transitive, reflexive cover of  $\prec\!\!\cdot$. 

\vspace{0.4cm}
\noindent \textbf{Comment 4. } ({\it $F$-partial ordered set})

\noindent Cobweb poset $\Pi = \langle\Phi,\leq\rangle$ is naturally   graded and sequence  $F$ - denominated thereby we call it sometimes  \textbf{$F$-partial ordered set  $\langle\Phi,\leq\rangle$}.

%%%%%%%%%%%%%%%%%%%%%%%
%%%%%%%%%%%%%%%%%%%%%%%%%%
%%%%%%%%%%%%%%%%%%%%%%%%

\section{ Dimension of cobweb posets-revisited. }

%%%%%%%%
\subsection{ oDAG \cite{7} }

\begin{observen} [cobwebs are oDAGs]
In  \cite{2}  it was observed that cobweb posets' Hasse diagrams   are the members of so called oDAGs family i.e.   cobweb posets' Hasse diagrams  are orderable Directed Acyclic  Graphs which is equivalent to say that the associated poset $P(D) = (\Phi, \leq)$ of   $D = (\Phi,\prec\!\!\cdot )$  of  is of dimension 2.
\end{observen}

\vspace{0.2cm}
\noindent \textbf{Recall:} DAGs - hence graded digraphs with minimal elements always  might be considered - up to digraphs isomorphism - as natural digraphs \cite{8} i.e. digraphs with natural   labeling (i.e.  $x_i < x_j  \Rightarrow  i < j$ ).

\begin{defn}[Plotnikov - see \cite{7} , \cite{2} and then below]
A digraph $D = (\Phi,\prec)$    is called the orderable digraph (oDAG) if there exists a dim 2 poset such  that its Hasse diagram coincides with the digraph $G$".
\end{defn}

\vspace{0.4cm}
\noindent The statement from \cite{2}  may be now restated as follows:

\begin{observen} [oDAG]
Cobweb  $P(D) =(\Phi , \leq)$   posets'  Hasse diagrams $D = (\Phi,\prec\!\!\cdot )$ are oDAGs.
\end{observen}

\vspace{0.2cm}
\noindent \emph{Proof}: Obvious. Cobweb posets are posets with minimal elements set $\Phi_0$. Cobweb posets Hasse diagrams are DAGs. Cobweb posets representing the natural join of  are then dim 2  posets as their Hasse digraphs are intersection of a natural labeling linear order $L_1$  and its "dual"  $L_2$ denominated correspondingly in a standard way by: 
$L_1$  = natural labeling:  chose for the topological ordering $L_1$ the labeling of minimal elements set $\Phi_0$  with  labels $1,2,...$,  from the left to the right  ( see Fig2. ) then proceed up to the next level  $\Phi_1$  and continue the labeling  "$\to$" from the left to the right  [$\Phi_1$ is now treated as the set o minimal elements if $\Phi_0$ is removed]  and so on. Apply the procedure of  subsequent removal of minimal elements  i.e. removal of subsequent labeled levels  $F_k$ - labeling the vertices along the levels  from the left to the right.  

\vspace{0.2cm}
\noindent $L_2$ = "dual" natural labeling:  chose for the topological ordering $L_2$ the labeling of minimal elements set $F_0$  with  labels $1,2,...$,  from the right  to  the left to ( see Fig1. ) then proceed up to the next level  $F_1$  and continue the labeling "$\leftarrow$" from the right to  the left  [$\Phi_1$   is now treated as the set o minimal elements if  $\Phi_0$    is removed]  and so on. Apply the procedure of  subsequent removal of minimal elements  i.e. removal of subsequent labeled levels   $\Phi_k$  - labeling now the vertices along the levels  from the right to  the left    q.e.d.

%%%%%%%%%%%%%
\subsection{ Brief history of the short  oDAG's name  life }

On the history of  oDAG nomenclature with David Halitsky and Others input one is expected to see  more in \cite{15}. See also the December $2008$
subject of  The Internet Gian Carlo Rota Polish Seminar ($http://ii.uwb.edu.pl/akk/sem/sem\_rota.htm$). Here we present its  sub-history  leading the author to note that cobweb posets are oDAGs.  

\vspace{0.2cm}
\noindent According to Anatoly Plotnikov the concept and the name  of oDAG was introduced  by David Halitsky from Cumulative Inquiry in 2004.

\vspace{0.2cm}

\vspace{0.2cm}

\noindent \textbf{oDAG-2004}  (Plotnikov)

\vspace{0.2cm}
\noindent Quote 1. \emph{"A digraph $G \in D_n$  will be called orderable  (oDAG) if there exists are dim 2 poset such  that its Hasse diagram coincide with the digraph $G$"}.  

\vspace{0.1cm}

\noindent The Quote 1 comes from \cite{9}  in  \cite{2}  i.e. from A.D. \emph{Plotnikov A formal approach to the oDAG/POSET problem} (2004) \emph{$html://www.cumulativeinquiry.com/Problems/solut2.pdf$}    (submitted to publication - March 2005)

\vspace{0.2cm}
\noindent The quote of the Quote 1 is to be found in  \cite{9}

\vspace{0.2cm}
\noindent oDAG-2005  \cite{2}

\vspace{0.1cm}

\noindent Quote 2  \emph{"A digraph G is called the orderable digraph (oDAG) if there exists a dim 2 poset such  that its Hasse diagram coincides with the digraph $G$"}. \cite{2}

\vspace{0.2cm}
\noindent oDAG-2006  \cite{7} 

\vspace{0.1cm}

\noindent Quote 3 \emph{"A digraph G is called the orderable  if there exists a dim 2 poset such  that its Hasse diagram coincides with the digraph  $G$"}.   \cite{7}

\vspace{0.2cm}
\noindent For further use of oDAG nomenclature see \cite{6}, and references therein. 
For further references and  recent  results on cobweb posets see  \cite{10} and  \cite{11}.

\begin{defn} [KoDAG]
The transitive and reflexive reduction of cobweb poset $\Pi = \langle\Phi,\leq\rangle$ i.e.   posets' $\Pi$ cover relation digraph [Hasse  diagram]  $D = (\Phi, \prec\!\!\cdot)$  is called KoDAG.   
\end{defn}

\noindent See \cite{11}-\cite{14}.

\vspace{0.2cm}
\noindent \textbf{Comment 5.} Apply Comment 1.  

\vspace{0.2cm}
\noindent \textbf{Why} do \textbf{we stick}  to call  KoDAGs graded digraphs with  associated poset  $\Pi = \langle\Phi,\leq\rangle$ the \textbf{orderable} DAGs  on their own independently of the nomenclature quoted ?

\vspace{0.2cm}
\noindent Let $D = (\Phi,\prec\!\!\cdot )$  denotes now  any transitive irreducible DAG [ for example any \textbf{graded} digraph including  KoDAG digraph for example as above].  Let poset  $P(D) = (\Phi, \leq)$  be associated to $D = (\Phi,\prec\!\!\cdot)$   

\begin{defn} [Ferrers dimension]
We say that the poset  $P(D) = (\Phi, \leq)$ is of Ferrers dimension $k$  iff  it is associated to $D = (\Phi,\prec\!\!\cdot )$  of   Ferrers dimension $k$.
\end{defn}

\begin{observen} [Ferrers dimension]
Cobweb posets are posets of Ferrers dimension equal to one.
\end{observen}

\noindent \emph{Proof.}  Apply any of many characterizations of Ferrers  digraphs to see that cobweb posets are posets' cover relation digraphs [Hasse diagrams] are  Ferrers   digraphs. For example consult Section 3 and see that biadjacency matrix does not contain  any of  two $ 2\times 2$ permutation matrices.

\vspace{0.4cm}
\noindent \textbf{Comment 6. }
Any KoDAG digraph $D = (\Phi,\prec\!\!\cdot )$ is the digraph stable under  the transitive and reflexive reduction i.e. ["`irreducible"'] Hasse portrait of  Ferrers relation $\prec\!\!\cdot$.  The positions of 1's in biadjacecy [reduced adjacency] matrix  display the support of Ferrers relation $\prec\!\!\cdot$.
$D = (\Phi,\prec\!\!\cdot)$ is then interval order relation digraph. The  digraph $(\Phi,\leq)$ of the cobweb poset  $P(D) = (\Phi, \leq)$    associated to KoDAG digraph $D = (\Phi,\leq )$ is the portrait of  Ferrers relation  $\leq$.   The positions of 1's in biadjacecy [reduced adjacency]  matrix display the support of Ferrers relation $\leq$. Note: for  $F$-denominated cobweb posets  the nomenclature identifies:  biajacency  [reduced adjacency] matrix  $\equiv$  zeta matrix i.e. the incidence matrix $\zeta_F$   of the $F$- poset (see: Fig.$\zeta_N$ and Fig.$\zeta_F$ ). Recall that this $F$-partial ordered set $\langle\Phi,\leq\rangle$   is a natural join of   $F$-chain of  binary $K$-relations (complete or universal relations as called sometimes). These relations are  represented by di-bicliques $\stackrel{\rightarrow}{K_{k,k+1}}$  which are on their own the   Ferrers dimension one digraphs. 
As for the other - not necessarily $K$-relations' chains we may end up with Ferrers or not digraphs in corresponding di-bigraphs'  chain.         
See below, then Section 4   and more  in  \cite{15}.

%%%%%%%%%%%%%%%%%%
\section{  The natural join $\os$ operation }

We define here the adjacency matrices representation of the natural join $\os$ operation.

%%%%%%%%%%%%
\subsection{ Recall }

Let  $D(R) = (V(R)\cup W(R),E(R)) \equiv (V \cup W, E) ;   V\cap W = \emptyset ,   E (R)  \subseteq V\times W$.
Let  $D(R)$ denotes here down the  \emph{bipartite digraph of binary relation} $R$   with  $\mathrm{dom}(R) = V$ and  $\mathrm{rang}(R)=W$.       Colligate with the  anticipated  examples  $R = R_k  \subseteq \Phi_k\times\Phi_{k+1}  \equiv \stackrel{\rightarrow}{K_{k,k+1}}, V(R)\cup W(R)= \Phi_k  \cup \Phi_{k+1}$.

%%%%%%%%%%%%%
\subsection{ The adjacency matrices and their  natural join.}

The adjacency matrix $\mathbf{A}$ of a bipartite graph with \textbf{biadjacency} (reduced adjacency \cite{16}) matrix $\mathbf{B}$ is given by

$$
	\mathbf{A} = \left(
	\begin{array}{cc}
	0 & \mathbf{B} \\
	\mathbf{B}^T & 0 \\
	\end{array}
	\right).
$$

\begin{defn}
The adjacency matrix $\mathbf{A}[D]$  of a bipartite digraph $D(R)= ( P\cup L ,  E \subseteq P\times L)$   with biadjacency matrix $\mathbf{B}$ is given by

$$
	\mathbf{A}[D] = \left(
	\begin{array}{cc}
	0_{k,k} & \mathbf{B}(k\times m) \\
	0_{m,k} & 0_{m,m} \\
	\end{array}
	\right).
$$
where  $k = | P |$,   $m = | L |$.
\end{defn}

\begin{conven}
 $S \copyright R$ = composition of  binary relations  $S$  and  $R  \Leftrightarrow   \mathbf{B}_{R\copyright S} =  \mathbf{B}_R \copyright \mathbf{B}_S$ where    ( $|V|= k , |W|= m$ )   $\mathbf{B}_R (k \times m) \equiv \mathbf{B}_R$  is the $(k \times m)$ 
\end{conven}

\noindent \textbf{biadjacency} [or another name:  \textbf{reduced}   adjacency]  matrix  of the bipartite relations' $R$  digraph $B(R)$  and $\copyright$ apart from  relations composition  denotes also  Boolean multiplication of these rectangular biadjacency  Boolean matrices  $B_R , B_S$.    What is their form?   The answer is in the block structure of  the  standard square $(n \times n)$ adjacency matrix $A[D(R)];  n = k +m$ .  The form of  standard square adjacency matrix $A[G(R)]$ of bipartite digraph  $D(R)$ has the following  apparently  recognizable block reduced structure:  [ $O_{s\times s}$ stays for $(k \times m)$ zero matrix ]

$$
	\mathbf{A}[D(R)] = \left[
	\begin{array}{ll}
		O_{k\times k} & \mathbf{A}_R(k\times m) \\
		O_{m\times k} & O_{m\times m}
	\end{array}
	\right]
$$

\noindent Let $D(S) = (W(S)\cup T(S),E(S))$; $W\cap T = \emptyset$, $E (S)  \subseteq W\times T;$ ($|W|= m, |T|= s$); hence

$$
	\mathbf{A}[D(S)] = \left[
	\begin{array}{ll}
		O_{m\times m} & \mathbf{A}_S(m\times s) \\
		O_{s\times m} & O_{s\times s}
	\end{array}
	\right]
$$

\begin{defn} [natural join condition]
The ordered pair of matrices   $\langle \mathbf{A_1}, \mathbf{A_2} \rangle$ is said to satisfy the natural join condition iff  they  have the block structure  of     $\mathbf{A}[D(R)]$  and  $\mathbf{A}[D(S)]$  as above  i.e. iff  they might be identified accordingly : $\mathbf{A_1} =  \mathbf{A}[D(R)]$   and  $\mathbf{A_2} =  \mathbf{A}[D(S)]$.   
\end{defn}

\noindent Correspondingly if  two given digraphs $G_1$ and $G_2$  are such that their adjacency matrices  $\mathbf{A_1} = \mathbf{A}[G_1]$  and  $\mathbf{A_2} =\mathbf{A}[G_2]$ do satisfy the natural join condition we shall say  that $G_1$ and $G_2$  satisfy the natural join condition.
 For matrices satisfying the natural join condition one may define what follows.

\vspace{0.4cm}
\noindent  First we  define the \textbf{Boolean reduced}  or \textbf{natural join  composition} $\cs$   and secondly the natural join $\os$ of adjacent matrices  satisfying the natural join condition. 

\begin{defn} ($\cs$ composition)

$$
	\mathbf{A}[D(R\copyright S)] =: \mathbf{A}[D(R)] \cs \mathbf{A}[D(S)] =  \left[
	\begin{array}{ll}
		O_{k\times k} & \mathbf{A}_{R\copyright S}(k\times s) \\
		O_{s\times k} & O_{s\times s}
	\end{array}
	\right]
$$

\noindent where $\mathbf{A}_{R\copyright S}(k\times s) = \mathbf{A}_R(k\times m) \copyright \mathbf{A}_S(m\times s)$.
\end{defn}

\noindent according  to the scheme: 
$$
	[(k+m) \times (k + m )]  \cs [(m + s) \times (m + s)]  =  [(k+ s) \times (k+ s)] .
$$

\vspace{0.4cm}
\noindent \textbf{Comment 7.}
\noindent The adequate projection makes out the intermediate, joint in  common  $\mathrm{dom}(S) = \mathrm{rang}(R)=W$ , $|W|= m$.

\vspace{0.4cm}
\noindent The above Boolean reduced composition $\cs$  of adjacent matrices technically reduces then to the calculation of just Boolean product of  the  \textbf{reduced}  rectangular  adjacency matrices  of the bipartite relations` graphs.

\vspace{0.2cm}
\noindent We  are however  now in need of the Boolean natural  join product  $\os$  of adjacent matrices  already announced at the beginning of this presentation. Let us now define it.

\vspace{0.4cm}
\noindent As for  the \textbf{natural join} notion we aim at the morphism correspondence:
$$
	S \os R  \Leftrightarrow    M_{S\os R}  =  M_R \os M_S
$$

\noindent where $S \os R$ = natural  join of  binary relations  $S$  and  $R$  while
$M_{S\os R} =  M_R \os M_S$ = natural  join of  standard square adjacency matrices  
(with   customary convention: $M[G(R)]  \equiv  M_R$  adapted). Attention:   recall here that  the natural join of the above binary relations  $R \os S$  is  the ternary relation - and on one results in $k$-ary relations  if with more factors undergo the $\os$ product.  As a matter of fact \textbf{ $\os$ operates on multi-ary relations according to the scheme:}   

$$
	(n+k)_{ary} \os (k+m)_{ary}  =   (n+ k +m)_{ary} .
$$

\noindent For example: $(1+1)_{ary} \os (1+1)_{ary} = (1+ 1 +1)_{ary}, binary \os binary = ternary$. 

\vspace{0.4cm}
\noindent Technically - the natural join of the $k$-ary  and  $n$-ary relations  is  defined accordingly the same way via  $\os$ natural join product of adjacency matrices - the adjacency matrices of   these relations' Hasse digraphs.

\vspace{0.2cm}
\noindent With the notation established above we finally define the natural join  $\os$  of two adjacency matrices as  follows:

\begin{defn} [natural join $\os$ of  biadjacency  matrices]. 

$$
 A[D(R \os S)]  =:  A[D(R)]  \os A[D(S)] =                                                
$$

$$
	= \left[
	\begin{array}{ll}
		O_{k\times k} & A_R(k\times m) \\
		O_{m\times k} & O_{m\times m}
	\end{array}
	\right] 
	\os
	\left[
	\begin{array}{ll}
		O_{m\times m} & A_S(m\times s) \\
		O_{s\times m} & O_{s\times s}
	\end{array}
	\right] =
$$
$$
	=\left[
	\begin{array}{lll}
		O_{k\times k} & A_R(k\times m) & O_{k\times s}\\
		O_{m\times k} & O_{m\times m}  & A_S(m\times s) \\
		O_{s\times k} & O_{s\times m}  & O_{s\times s}
	\end{array} 
	\right]
$$
\end{defn}

\vspace{0.4cm}
\noindent \textbf{Comment 8}.  The adequate projection used in natural join operation lefts one copy of the joint in  common "intermediate" submatrix $O_{m\times m}$ and consequently lefts one copy  of  "intermediate" joint in  common  $m$ according  to the scheme:
$$
	[(k+m) \times (k + m )]  \os  [(m + s) \times (m + s)]  =  [(k+ m + s) \times (k+ m + s)]  .
$$

%%%%%%%%%%
\subsection{  The biadjacency matrices of  the  natural join of adjacency matrices.  }

Denote  with  $B(A)$ the   biadjacency matrix of the adjacency matrix $A$.

\vspace{0.2cm}
\noindent Let  $A(G)$ denotes the adjacency matrix of the digraph $G$  , for example a di-biclique relation digraph.   Let  $A(G_k)$, $k= 0,1,2,...$ be the sequence adjacency matrices of the sequence $G_k, k=0,1,2,...$ of digraphs.  Let us identify  $B(A)\equiv B(G)$ as a convention.

\begin{defn} [digraphs natural join]
Let  digraphs $G_1$ and $G_2$   satisfy the natural join condition.  Let us make then  the identification  $A(G_1 \os G_2)  \equiv A_1 \os A_2$  as definition.  The digraph  $G_1 \os G_2$  is called the digraphs natural join of  digraphs  $G_1$ and $G_2$. Note that the order is essential. 
\end{defn}

\vspace{0.2cm}
\noindent We observe at once what follows. 

\begin{observen}
$$
	B (G_1 \os G_2) \equiv B (A_1 \os A_2) =  B(A_1)\oplus B(A_2) \equiv B (G_1)\oplus B(G_2)
$$
\end{observen}

\vspace{0.2cm}
\noindent \textbf{Comment 9.} The Observation 4 justifies the notation  $\os$ for the natural join
of relations digraphs and  equivalently for  the natural join of their adjacency matrices and 
equivalently for the natural join of relations that these are faithful representatives of.

\vspace{0.2cm}
\noindent As a consequence we have.

\begin{observen}
$$
	B\left(\os_{i=1}^n\right) \equiv  B [\os_{i=1}^n  A(G_i)]  =  \oplus_{i=1}^n  B[A(G_i) ]  \equiv  \mathrm{diag} (B_1 , B_2 , ..., B_n) = 
$$
$$
	= \left[ \begin{array}{lllll}
	B_1 \\
	& B_2 \\
	& & B_3 \\
	& ... & ... & ...\\
	& & & & B_n
	\end{array} \right]
$$

\noindent $n \in N \cup \{\infty\}$.  
\end{observen}

%%%%%%%%%%%%%%%
\subsection{  Applications }

Once any natural number valued sequence $F = \{F_n\}_{n\geq 1}$  is being chosen its KoDAG digraph is identified with Hasse cover relation  digraph.  Its adjacency matrix  $\mathbf{A}_F$   is sometimes called Hasse matrix  and is given in a plausible form and impressively straightforward  way. Just use the fact that  the Hasse digraph which is displaying cover relation  $\prec\!\!\cdot$  is  an  $F$  -chain of coined  bipartite digraphs -  coined each preceding with a subsequent one by  natural join operator  $\os$ [resemblance of  $\os$  to  direct matrix  sum is not naive - compare "natural join" of disjoint digraphs   with no common set of  marked nodes ("attributes") ].

\vspace{0.1cm}
\noindent Note: $I (s\times k)$  stays for  $(s\times k)$  matrix  of  ones  i.e.  $[ I (s\times k) ]_{ij} = 1$;  $1 \leq i \leq s,  1\leq j  \leq k$. 

\vspace{0.2cm}
\noindent Let us start  first with $F = \{F_n\}_{n\geq 1} = N$.  See \textbf{Fig.2} . Then its associated  \textbf{$F$-partial ordered set $\langle\Phi,\leq\rangle$} has the following Hasse digraph displaying  cover relation of the $\leq$  partial order

%%%%%%% picture 
\begin{figure}[ht]
\begin{center}
	\includegraphics[width=80mm]{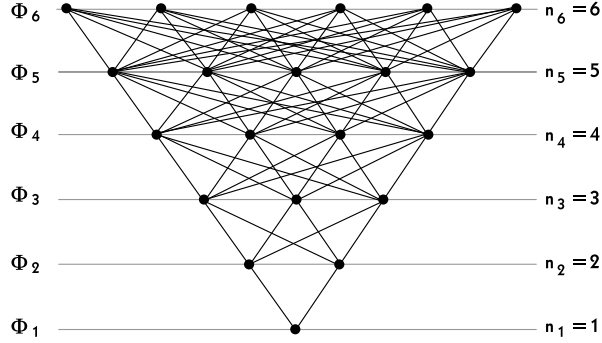}
	\caption {Display  of  a finite subposet $\Pi_6$ of the $N$ natural numbers cobweb poset}
\end{center}
\end{figure}

\vspace{0.2cm}
\noindent The Hasse matrix $\mathbf{A}_N$ i.e. adjacency matrix of cover relation digraph i.e. adjacency matrix  of the Hasse diagram of  the $N$-denominated cobweb poset $\langle\Phi,\leq\rangle$ is given by  upper triangular matrix $\mathbf{A}_N$ of the form:

$$
	\mathbf{A}_N = 
	\left[ \begin{array}{llllll}
	O_{1\times 1} & I(1\times 2)  & O_{1\times \infty} \\
	O_{2\times 1} & O_{2\times 2} & I(2\times 3)  & O_{2 \times \infty} \\
	O_{3\times 1} & O_{3\times 2} & O_{3\times 3} & I(3\times 4)  & O_{3 \times \infty} \\
	O_{4\times 1} & O_{4\times 2} & O_{4\times 3} & O_{4\times 4} & I(4\times 5) & O_{4 \times \infty} \\
	... etc & ... & and & so & on ...
	\end{array} \right]
$$

\noindent One may see that the zeta  function matrix of the $F = N$  choice is geometrical series in  $\mathbf{A}_N$  i.e. the geometrical series
in  the poset $\langle\Phi,\leq\rangle$ Hasse matrix $\mathbf{A}_N$:

$$
	\zeta = (1 - \mathbf{A}_N)^{-1 \copyright}
$$

\vspace{0.2cm}
\noindent Explicitly: $\zeta = (1-\mathbf{A}_N)^{-1\copyright} \equiv I_{\infty \time \infty} + \mathbf{A}_N + \mathbf{A}_N ^{\copyright 2} + ... = $

$$
	= \left[ \begin{array}{lllll}
	I_{1\times 1} & I(1\times \infty) \\
	O_{2\times 1} & I_{2\times 2} & I(2\times\infty) \\
	O_{3\times 1} & O_{3\times 2} & I_{3\times 3} & I(3\times\infty) \\
	O_{4\times 1} & O_{4\times 2} & O_{4\times 3} & I_{4\times 4} & I(4\times\infty) \\
	... etc & ... & and & so & on ...
	\end{array}\right]
$$

\vspace{0.2cm}
\noindent $\zeta = (1-\mathbf{A}_N)^{-1\copyright}$ because [let $\mathbf{A}_N = \mathbf{A}$]

\vspace{0.2cm}
\noindent $\mathbf{A}^k_{ij} =$ the number of maximal $k$-chains [$k>0$]  from the $x_0 \in \Phi_i$ to  $x_k \in \Phi_j$ i.e. here

$$
	\mathbf{A}^k_{ij} = \left\{ 
	\begin{array}{ll}
	0 & k\neq j-i \\
	\frac{j!}{i!} & k=j-k
	\end{array}
 	\right.  
	\mathrm{ hence }\ \ 
	\mathbf{A}^{\copyright k}_{ij} = \left\{ 
	\begin{array}{ll}
	1 & k = j-i \\
	0 & k \neq j-k
	\end{array}
 	\right. .
$$

\noindent and the supports (\emph{nonzero matrices blocks}) of  $\mathbf{A}^{\copyright k}$  and  $\mathbf{A}^{\copyright m}$  are disjoint  for $k \neq m$. Indeed:   the entry in row $i$ and column $j$ of  the inverse $(I - \mathbf{A})^{-1}$ gives \emph{the number of directed paths} from vertex $x_i$ to vertex $x_j$.                                                                                       This can be seen from geometric series with adjacency matrix as an argument 

$$
	(I - \mathbf{A})^{-1} = I + \mathbf{A} + \mathbf{A}^2 + \mathbf{A}^3 + ...
$$

\noindent taking care of  the fact that the number of paths from $i$ to $j$ equals the number of paths of length $0$ plus the number of paths of length $1$ plus the number of paths of length $2$, etc.

\vspace{0.2cm}
\noindent Therefore the entry in row $i$ and column $j$ of  the inverse $(I - \mathbf{A})^{-1\copyright}$ gives  the answer whether there exists a \emph{directed paths} from vertex $i$ to vertex $j$  (Boolean value 1) or not (Boolean value 0)  i.e. whether these vertices are comparable i.e.  whether $x_i < x_j$ or not.

\vspace{0.2cm}
\noindent \textbf{Remark:} 
In the cases - Boolean  poset $2^N$ and the "Ferrand-Zeckendorf"  poset of finite 
subsets of $N$  without two consecutive elements considered in \cite{17} one has

$$
	\zeta = exp[\mathbf{A}] = (1-\mathbf{A})^{-1\copyright} \equiv I_{\infty\times\infty} + \mathbf{A} + \mathbf{A}^{\copyright 2} + ...
$$

\noindent because in those cases

$$
	\mathbf{A}^k_{ij} = \left\{ 
	\begin{array}{ll}
	0  & k\neq j-i \\
	k! & k = j - k
	\end{array}
 	\right.  
	\mathrm{ hence }\ \ 
	\frac{1}{k!} \mathbf{A}^k_{ij} = \mathbf{A}^{\copyright k}_{ij} = \left\{ 
	\begin{array}{ll}
	1 & k = j-i \\
	0 & k \neq j-k
	\end{array}
 	\right. .
$$

\noindent How it goes in our $F$-case? Just see $\mathbf{A}_N^{\copyright 2}$ and  then  add $\mathbf{A}_N^{\copyright 0} \vee \mathbf{A}_N^{\copyright 1} \vee \mathbf{A}_N^{\copyright 2} \vee ...$

\noindent For example:

$$
	\mathbf{A}_N^{\copyright 2} = \left[ \begin{array}{lllllll}
	O_{1\times 1} & O_{1\times 2}  & I(1\times 3)  & O_{1\times \infty} \\
	O_{2\times 1} & O_{2\times 2}  & O_{2\times 3} & I(2 \times 4) & O_{2\times \infty} \\
	O_{3\times 1} & O_{3\times 2}  & O_{3\times 3} & O_{3\times 4} & I(3 \times 5) & O_{3\times \infty} \\
	O_{4\times 1} & O_{4\times 2}  & O_{4\times 3} & O_{4\times 4} & O_{4\times 5} & I(4 \times 6) & O_{4\times\infty} \\
... etc & ... & and & so & on ...
	\end{array}\right]
$$

\noindent Consequently  we arrive at  the incidence matrix $\zeta = \mathrm{exp}[\mathbf{A}_N]$ for the natural numbers cobweb poset displayed by Fig 3.  Note that incidence matrix $\zeta$ representing uniquely its corresponding  cobweb poset does  exhibits (see below) a staircase structure of zeros above the diagonal which is characteristic to Hasse diagrams of \textbf{all} cobweb posets. 

\vspace{1mm}
$$ \left[\begin{array}{ccccccccccccccccc}
1 & 1 & 1 & 1 & 1 & 1 & 1 & 1 & 1 & 1 & 1 & 1 & 1 & 1 & 1 & 1 & \cdots\\
0 & 1 & 0 & 1 & 1 & 1 & 1 & 1 & 1 & 1 & 1 & 1 & 1 & 1 & 1 & 1 & \cdots\\
0 & 0 & 1 & 1 & 1 & 1 & 1 & 1 & 1 & 1 & 1 & 1 & 1 & 1 & 1 & 1 & \cdots\\
0 & 0 & 0 & 1 & 0 & 0 & 1 & 1 & 1 & 1 & 1 & 1 & 1 & 1 & 1 & 1 & \cdots\\
0 & 0 & 0 & 0 & 1 & 0 & 1 & 1 & 1 & 1 & 1 & 1 & 1 & 1 & 1 & 1 & \cdots\\
0 & 0 & 0 & 0 & 0 & 1 & 0 & 0 & 1 & 1 & 1 & 1 & 1 & 1 & 1 & 1 & \cdots\\
0 & 0 & 0 & 0 & 0 & 0 & 1 & 0 & 0 & 0 & 1 & 1 & 1 & 1 & 1 & 1 & \cdots\\
0 & 0 & 0 & 0 & 0 & 0 & 0 & 1 & 0 & 0 & 1 & 1 & 1 & 1 & 1 & 1 & \cdots\\
0 & 0 & 0 & 0 & 0 & 0 & 0 & 0 & 1 & 0 & 0 & 0 & 0 & 1 & 1 & 1 & \cdots\\
0 & 0 & 0 & 0 & 0 & 0 & 0 & 0 & 0 & 1 & 1 & 0 & 0 & 1 & 1 & 1 & \cdots\\
0 & 0 & 0 & 0 & 0 & 0 & 0 & 0 & 0 & 0 & 1 & 0 & 0 & 0 & 0 & 1 & \cdots\\
0 & 0 & 0 & 0 & 0 & 0 & 0 & 0 & 0 & 0 & 0 & 1 & 0 & 0 & 0 & 1 & \cdots\\
0 & 0 & 0 & 0 & 0 & 0 & 0 & 0 & 0 & 0 & 0 & 0 & 1 & 0 & 0 & 1 & \cdots\\
0 & 0 & 0 & 0 & 0 & 0 & 0 & 0 & 0 & 0 & 0 & 0 & 0 & 1 & 0 & 1 & \cdots\\
0 & 0 & 0 & 0 & 0 & 0 & 0 & 0 & 0 & 0 & 0 & 0 & 0 & 0 & 1 & 1 & \cdots\\
0 & 0 & 0 & 0 & 0 & 0 & 0 & 0 & 0 & 0 & 0 & 0 & 0 & 0 & 0 & 1 & \cdots\\
. & . & . & . & . & . & . & . & . & . & . & . & . & . & . & . & . \cdots\\
 \end{array}\right]$$
\vspace{1mm}   \noindent \textbf{Figure $\zeta_N$.  The incidence matrix
$\zeta$ for the  natural numbers  i.e. N- cobweb poset}

\vspace{0.2cm}
\noindent \textbf{Comment  9.} 
 The given  $F$-denominated staircase zeros structure above the diagonal of zeta matrix $zeta$ is the \textbf{unique characteristics} of  its corresponding  \textbf{$F$-KoDAG} Hasse digraphs.

\vspace{0.2cm}
\noindent For example see Fig $\zeta_F$. below (from \cite{6}).

\vspace{1mm}
$$ \left[\begin{array}{ccccccccccccccccc}
1 & 1 & 1 & 1 & 1 & 1 & 1 & 1 & 1 & 1 & 1 & 1 & 1 & 1 & 1 & 1 & \cdots\\
0 & 1 & 1 & 1 & 1 & 1 & 1 & 1 & 1 & 1 & 1 & 1 & 1 & 1 & 1 & 1 & \cdots\\
0 & 0 & 1 & 1 & 1 & 1 & 1 & 1 & 1 & 1 & 1 & 1 & 1 & 1 & 1 & 1 & \cdots\\
0 & 0 & 0 & 1 & 0 & 1 & 1 & 1 & 1 & 1 & 1 & 1 & 1 & 1 & 1 & 1 & \cdots\\
0 & 0 & 0 & 0 & 1 & 1 & 1 & 1 & 1 & 1 & 1 & 1 & 1 & 1 & 1 & 1 & \cdots\\
0 & 0 & 0 & 0 & 0 & 1 & 0 & 0 & 1 & 1 & 1 & 1 & 1 & 1 & 1 & 1 & \cdots\\
0 & 0 & 0 & 0 & 0 & 0 & 1 & 0 & 1 & 1 & 1 & 1 & 1 & 1 & 1 & 1 & \cdots\\
0 & 0 & 0 & 0 & 0 & 0 & 0 & 1 & 1 & 1 & 1 & 1 & 1 & 1 & 1 & 1 & \cdots\\
0 & 0 & 0 & 0 & 0 & 0 & 0 & 0 & 1 & 0 & 0 & 0 & 0 & 1 & 1 & 1 & \cdots\\
0 & 0 & 0 & 0 & 0 & 0 & 0 & 0 & 0 & 1 & 0 & 0 & 0 & 1 & 1 & 1 & \cdots\\
0 & 0 & 0 & 0 & 0 & 0 & 0 & 0 & 0 & 0 & 1 & 0 & 0 & 0 & 1 & 1 & \cdots\\
0 & 0 & 0 & 0 & 0 & 0 & 0 & 0 & 0 & 0 & 0 & 1 & 0 & 1 & 1 & 1 & \cdots\\
0 & 0 & 0 & 0 & 0 & 0 & 0 & 0 & 0 & 0 & 0 & 0 & 1 & 1 & 1 & 1 & \cdots\\
0 & 0 & 0 & 0 & 0 & 0 & 0 & 0 & 0 & 0 & 0 & 0 & 0 & 1 & 0 & 0 & \cdots\\
0 & 0 & 0 & 0 & 0 & 0 & 0 & 0 & 0 & 0 & 0 & 0 & 0 & 0 & 1 & 0 & \cdots\\
0 & 0 & 0 & 0 & 0 & 0 & 0 & 0 & 0 & 0 & 0 & 0 & 0 & 0 & 0 & 1 & \cdots\\
. & . & . & . & . & . & . & . & . & . & . & . & . & . & . & . & . \cdots\\
 \end{array}\right]$$

\vspace{1mm} \noindent \textbf{Figure $\zeta_F$.  The incidence matrix
$\zeta$ for the Fibonacci cobweb poset associated to \textbf{$F$-KoDAG} Hasse digraph }

\vspace{0.2cm}
\noindent The  zeta matrix i.e. the incidence matrix $\zeta_F$  for the Fibonacci numbers cobweb poset \textbf{[$F$ - KoDAG]}  determines completely its incidence algebra  and corresponds to the poset  with Hasse diagram displayed by the Fig. 3.

%%%%%%% picture 
\begin{figure}[ht]
\begin{center}
	\includegraphics[width=80mm]{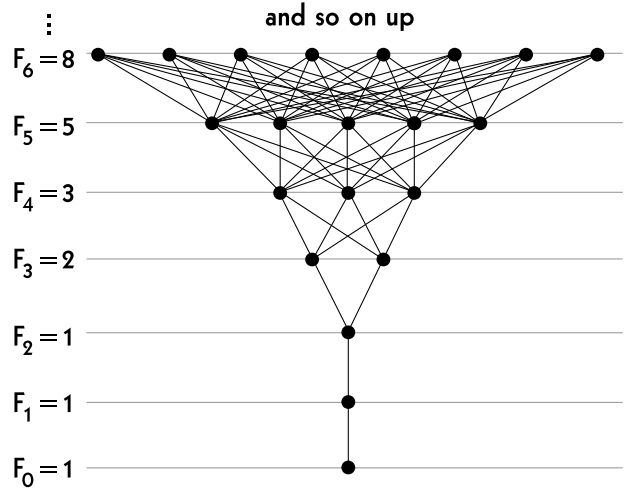}
	\caption {Display  of  the $F$- Fibonacci  numbers cobweb poset}
\end{center}
\end{figure}

\vspace{0.2cm}
\noindent The explicit  expression for zeta matrix $\zeta_F$  via known blocks of zeros and ones for arbitrary natural numbers valued $F$- sequence  is readily found due to brilliant mnemonic  efficiency  of the authors up-side-down notation (see Appendix in \cite{13}).  With this notation inspired by Gauss  and the reasoning just repeated   with   "$k_F$"  numbers   replacing  $k$ - natural numbers one gets in the spirit of Knuth \cite{18} the clean result:

$$
	\mathbf{A}_F = \left[\begin{array}{llllll}
	0_{1_F\times 1_F} & I(1_F \times 2_F) & 0_{1_F \times \infty} \\
	0_{2_F\times 1_F} & 0_{2_F\times 2_F} & I(2_F \times 3_F) & 0_{2_F \times \infty} \\
	0_{3_F\times 1_F} & 0_{3_F\times 2_F} & 0_{3_F\times 3_F} & I(3_F \times 4_F) & 0_{3_F \times \infty} \\
	0_{4_F\times 1_F} & 0_{4_F\times 2_F} & 0_{4_F\times 3_F} & 0_{4_F\times 4_F} & I(4_F \times 5_F) & 0_{4_F \times \infty} \\
	... & etc & ... & and\ so\ on & ...
	\end{array}\right]
$$

\noindent and

$$
	\zeta_F = exp_\copyright[\mathbf{A}_F] \equiv (1 - \mathbf{A}_F)^{-1\copyright} \equiv I_{\infty\times\infty} + \mathbf{A}_F + \mathbf{A}_F^{\copyright 2} + ... =
$$
$$
	= \left[\begin{array}{lllll}
	I_{1_F\times 1_F} & I(1_F\times\infty) \\
	O_{2_F\times 1_F} & I_{2_F\times 2_F} & I(2_F\times\infty) \\
	O_{3_F\times 1_F} & O_{3_F\times 2_F} & I_{3_F\times 3_F} & I(3_F\times\infty) \\
	O_{4_F\times 1_F} & O_{4_F\times 2_F} & O_{4_F\times 3_F} & I_{4_F\times 4_F} & I(4_F\times\infty) \\
	... & etc & ... & and\ so\ on & ...
	\end{array}\right]
$$

\vspace{0.2cm}
\noindent \textbf{Comment 10.} (ad "upside down notation")

\noindent Concerning Gauss and Knuth - see remarks in \cite{18} on Gaussian  binomial coefficients.

\begin{observen}
Let us denote by  $\langle\Phi_k\to\Phi_{k+1}\rangle$ (see the authors papers quoted) the di-bicliques  denominated by subsequent levels $\Phi_k, \Phi_{k+1}$ of the graded  $F$-poset $P(D) = (\Phi, \leq)$  i.e. levels $\Phi_k , \Phi_{k+1}$ of  its cover relation graded digraph  $D = (\Phi,\prec\!\!\cdot$)  [Hasse diagram].   Then

$$
	B\left(\os_{k=1}^n \langle\Phi_k\to\Phi_{k+1}\rangle \right) = \mathrm{diag}(I_1,I_2,...,I_n) = 
$$
$$
	= \left[ \begin{array}{lllll}
	I(1_F\times 2_F) \\
	& I(2_F\times 3_F) \\
	& & I(3_F\times 4_F) \\
	& & ... \\
	& & & & I(n_F\	times (n+1)_F)
	\end{array} \right]
$$

\noindent where $I_k \equiv I(k_F \times (k+1)_F)$, $k = 1,...,n$ and where  - recall - $I (s\times k)$  stays for $(s\times k)$  matrix  of  ones  i.e.  $[ I (s\times k) ]_{ij} = 1$;  $1 \leq i \leq  s,  1\leq j  \leq k.$  and  $n \in N \cup \{\infty\}$.  
\end{observen}

\begin{observen}
Consider bigraphs'  chain obtained from the above di-biqliqes' chain via  deleting or no  arcs making thus [if deleting arcs] some or all of the di-bicliques $ \langle\Phi_k\to\Phi_{k+1}\rangle$  not di-biqliques; denote  them as  $G_k$. Let $B_k = B(G_k)$ denotes their biadjacency matrices correspondingly.  Then for any such $F$-denominated chain [hence any chain ] of bipartite digraphs  $G_k$  the general formula is:

$$
 B\left( \os_{i=1}^n G_i \right) \equiv  B [\os_{i=1}^n  A(G_i)] =  \oplus_{i=1}^n  B[A(G_i) ]  \equiv  \mathrm{diag} (B_1 , B_2 , ..., B_n) =
$$
$$
	= \left[ \begin{array}{lllll}
	B_1 \\
	& B_2 \\
	& & B_3 \\
	& & ... \\
	& & & & B_n
	\end{array} \right]
$$

\noindent $n \in N \cup \{\infty\}$.
\end{observen}

\begin{observen}
 The $F$-poset $P(G) = (\Phi, \leq)$  i.e.  its cover relation graded digraph $G = (\Phi,\prec\!\!\cdot) = \os_{k=0}^m G_k$ is of Ferrers dimension one iff  in the process of deleting arcs from  the cobweb poset Hasse diagram  $D = (\Phi,\prec\!\!\cdot)$ = $\os_{k=0}^n  \langle\Phi_k\to\Phi_{k+1}\rangle $   does not produces  $2\times 2$  permutation submatrices in any  bigraphs $G_k$  biadjacency  matrix $B_k= B (G_k)$.
\end{observen}

\vspace{0.2cm}
\noindent \textbf{Examples} (finite subposets of cobweb posets) 

\noindent Fig.4 and Fig.5  display a Hasse diagram portraits  of finite subposets of cobweb posets. In view of the \textbf{Observation 2}  these subposets are naturally Ferrers digraphs i.e. of Ferrers dimension equal one.

%%%%%%% picture 
\begin{figure}[ht]
\begin{center}
	\includegraphics[width=100mm]{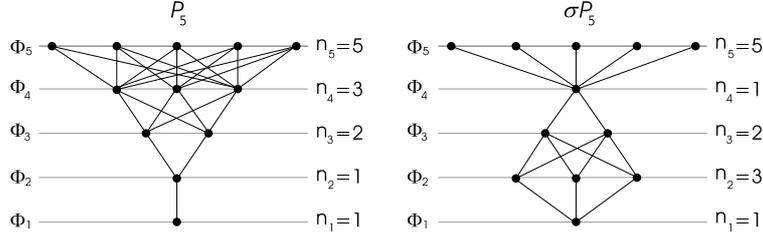}
	\caption {Display of the subposet $P_5$  of the $F$= Fibonacci sequence $F$-cobweb poset and $\sigma P_5$ subposet of the $\sigma$ permuted Fibonacci $F$-cobweb poset}
\end{center}
\end{figure}

%%%%%%% picture 
\begin{figure}[ht]
\begin{center}
	\includegraphics[width=100mm]{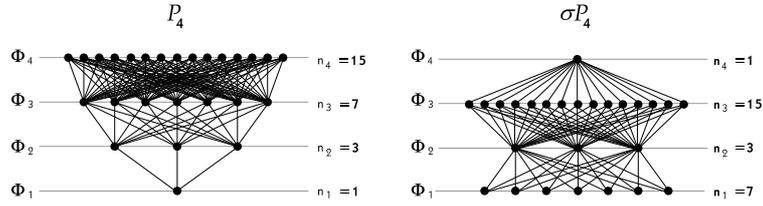}
	\caption {Display of the subposet $P_4$ of the  $F$ = Gaussian integers sequence $(q=2)$ $F$-cobweb poset and $\sigma P_4$ subposet of the $\sigma$ permuted Gaussian $(q=2)$ $F$-cobweb poset.}
\end{center}
\end{figure}

%%%%%%%%%%%%%%%%%%%%%%%%%%%%%%
\section{Summary}

\subsection{Principal - natural  identifications}

Any \textbf{KoDAG}  is  a \textbf{di}-bicliques chain  $\Leftrightarrow$  Any  \textbf{KoDAG is a natural join} of complete bipartite \textbf{graphs} [ \textbf{di}-bicliques ]  = 
$$
	( \Phi_0 \cup \Phi_1 \cup ... \cup \Phi_n \cup ..., E_0 \cup E_1\cup ... \cup E_n \cup ...) \equiv D(\bigcup_{k\geq 0}\Phi_k,\bigcup_{k\geq 0} E_k )  \equiv D (\Phi,E) 
$$

\noindent where $E_k = \Phi_k\times \Phi_{k+1} \equiv \stackrel{\rightarrow}{K_{k,k+1}}$ and $E = \bigcup_{k\geq 0}E_k$.

\vspace{0.2cm}
\noindent Naturally, as indicated earlier any graded posets' Hasse diagram  with finite width  including  \textbf{KoDAGs}  is of the form  
$$
	D (\Phi , E)  \equiv D( \bigcup_{k\geq 0}\Phi_k,\bigcup_{k\geq 0} E_k ) \Leftrightarrow  \langle \Phi,\leq \rangle
$$

\noindent where  $E_k \subseteq \Phi_k\times \Phi_{k+1} \equiv  \stackrel{\rightarrow}{K_{k,k+1}}$ and the definition of  $\leq$  from  \textbf{1.3.} is applied.
\noindent In front of all the above presentation the following is clear . 

\vspace{0.1cm}
\begin{observen}
"Many"   graded digraphs with finite width including \textbf{KoDAGs} $D = (V,\prec\!\!\cdot )$ encode bijectively  their  correspondent $n$-ary relation  ($n \in N \cup \{\infty\}$ as seen from its  following  definition: $ E_k \subseteq \Phi_k \times \Phi_{k+1} \equiv  \stackrel{\rightarrow}{K_{k,k+1}}$  where \\
\textcolor{red}{\textbf{($n$-ary relation)}}  $E = \os_{k=0}^{n-1} E_k \subset \times_{k=0}^n \Phi_k$ \\ 
i.e. identified with graded poset $\left\langle V_n, E \right\rangle$  natural join obtained $n+1$-ary relation $E$ is a subset of Cartesian product obtained the universal $n+1$-ary relation identified with cobweb poset digraph $\left\langle V_n,\prec\!\!\cdot \right\rangle$). $V_{\infty}\equiv V$.
\end{observen}

\vspace{0.2cm}
\noindent Which are those "many"? The characterization is arrived at   with au rebour point of view.
Any $n$-ary relation  ($n \in N \cup \{\infty\}$)  determines uniquely [may be identified with] its correspondent graded digraph with minimal elements set $\Phi_0$ given by the \textcolor{red}{\textbf{($n$-ary rel.)}}  formula  
$$E = \os_{k=0}^{n-1} E_k \subset \times_{k=0}^n \Phi_k,$$
where the sequence of binary relations $E_k \subseteq \Phi_k\times \Phi_{k+1} \equiv \stackrel{\rightarrow}{K_{k,k+1}}$ is denominated by the source $n$-ary relation as the following example shows.

\vspace{0.4cm}
\noindent \textbf{Example} (ternary = $Binary_1$ $\os$ $Binary_2$)

\noindent Let  $T \subset X\times Z\times Y$  where  $X =\{ x_1,x_2,x_3\}$, $Z = \{ z_1,z_2,z_3,z_4\}$, $Y = \{y_1,y_2\}$  and  
$$
	T = \{ \langle x_1,z_1,y_1 \rangle, \langle x_1,z_2,y_1 \rangle,   \langle x_1,z_4,y_2 \rangle, \langle x_2,z_3,y_2 \rangle, \langle x_3,z_3,y_2 \rangle \}.
$$

%%%%%%% picture 
\begin{figure}[ht]
\begin{center}
	\includegraphics[width=50mm]{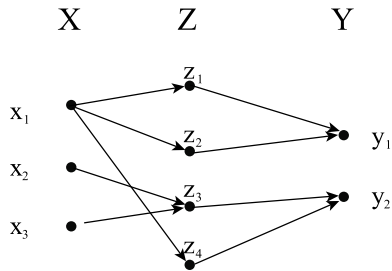}
	\caption {Display of example ternary = $Binary_1$ $\os$ $Binary_2$. \label{fig:ternary}}
\end{center}
\end{figure}

\noindent Let $X\times Z \supset E_1= \{ \langle x_1,z_1 \rangle, \langle x_1,z_2 \rangle, \langle x_1,z_4 \rangle, \langle x_2,z_3 \rangle, \langle x_3,z_3 \rangle \}$  and 
$Z\times Y \supset E_2 = \{ \langle z_1,y_1 \rangle, \langle z_2,y_1\rangle,  \langle z_3,y_1\rangle, \langle z_4,y_2\rangle \}$.
Then $T = E_1 \os E_2$.

\vspace{0.2cm}
\noindent More on that - see \cite{15} and  see references to the authors recent papers therein. 

\vspace{0.4cm}
\noindent \textbf{Comment 11.}  As a comment to the \textbf{Observation 9} and the \textbf{Observation 3} consider Fig.7  which was the source of inspiration for cobweb posets birth \cite{4,3,2,5,6} and here serves as Hasse diagram $D_{Fib} \equiv (\Phi, \prec\!\!\cdot_{Fib})$  of the poset $P(D_{Fib}) = (\Phi, \leq_{Fib} )$  associated to  $D_{Fib}$. Obviuosly, $P(D_{Fib})$  is a subposet of the Fibonacci cobweb poset $P(D)$ and  $D_{Fib}$  is a subgraph of the Fibonacci cobweb poset $P(D)$ Hasse diagram $D \equiv  (\Phi,\prec\!\!\cdot )$. 

\noindent The Ferrers dimension of $D_{Fib}$ is obviously not equal one. 

%%%%%%% picture 
\begin{figure}[ht]
\begin{center}
	\includegraphics[width=60mm]{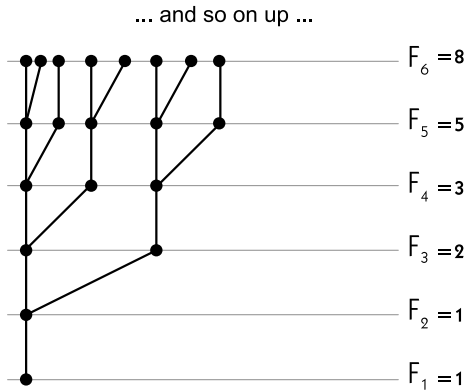}
	\caption {Display  of  of Hasse diagram of the form of the Fibonacci tree.}
\end{center}
\end{figure}

\vspace{0.2cm}
\noindent \textbf{Exercise.} Find the Ferrers dimension of  $D_{Fib}$. What is the dimension of the poset $P(D_{Fib}) = (\Phi, \leq_{Fib})$ ?   (Compare with \textbf{Observation 2}). Find the chain  $E_k \subset \Phi_{k}\times\Phi_{k+1}$, $k =0,1,2,...$ of binary relations such  that $D_{Fib,n} = \os_{k=0}^n E_k, n \in N \cup \{\infty\}$. Find the Ferrers dimension of $D_{Fib,n}$.

\vspace{0.4cm}
\noindent \textbf{Ad Bibliography Remark}

\noindent On the history of  oDAG nomenclature with David Halitsky and Others input one is expected to see  more in \cite{15}. See also the December $2008$ subject of  The Internet Gian Carlo Rota Polish Seminar ($http://ii.uwb.edu.pl/akk/sem/sem\_rota.htm$).
Recommended readings on Ferrers  digraphs of immediate use here are
\cite{19}-\cite{25}. For example see pages 61 an 85 in \cite{19}, see page 2 in \cite{20}. The  J. Riguet paper  \cite{21} is the source paper including also equivalent characterizations of Ferrers digraphs as well as other \cite{22,23,24}. The now classic reference on interval orders and interval graphs is \cite{25}.

\vspace{0.4cm}
\noindent \textbf{Acknowledgments} Thank are expressed here to the Student of Gda\'nsk University Maciej Dziemia\'nczuk for applying his skillful   TeX-nology with respect to the present work as well as for his general assistance and cooperation on KoDAGs  investigation.

\end{document}